\documentclass[11pt]{article}

\usepackage{graphicx}          

\usepackage{epsfig}
\usepackage{latexsym}
\usepackage{pdfsync}
\usepackage{amssymb}

\def\Real{{I\!\!R}}

\newcommand{\eq}{\begin{equation}\begin{array}{rllllllllllllllllllllllllllllllll}}
\newcommand{\ee}{\end{array}\end{equation}}
\newcommand{\bmt}{\left[ \begin{array}{ccccccccc}}
\newcommand{\emt}{\end{array}\right]}
\newcommand{\bea}{\begin{eqnarray}}
\newcommand{\eea}{\end{eqnarray}}
\newcommand{\bean}{\begin{eqnarray*}}
\newcommand{\eean}{\end{eqnarray*}}

\title{Adaptive Horizon Model Predictive Control\\ and Al'brekht's Method} 

\author{Arthur J Krener\\
Department of Applied Mathematics,\\ Naval Postgraduate School\\Monterey, CA 93943, USA\\ (email: ajkrener@nps.edu).}

\begin{document}
\maketitle

\begin{abstract}
A standard way of finding a feedback law that stabilizes a control system to an operating point
is to recast the problem as an infinite horizon optimal control problem.  If the optimal cost and the optmal feedback 
can be found on a large domain around the operating point then a Lyapunov argument can be used to verify the asymptotic stability
of the closed loop dynamics.  The problem with this approach is that is usually very difficult to find 
the optimal cost and the optmal feedback on a  large domain for nonlinear problems with or without constraints.
Hence the increasing interest in Model Predictive Control (MPC).  In standard MPC a  finite horizon optimal control problem  
is solved in real time but just at the current state, the first control action is implimented, the system evolves one time step and the process is repeated. 
  A terminal cost and terminal feedback found by Al'brekht's methoddefined in a neighborhood of the operating point is used to shorten the horizon and thereby make the nonlinear  programs easier to solve because they have less decision variables.    

Adaptive Horizon Model Predictive Control  (AHMPC) is a scheme for varying  the horizon length of Model Predictive Control (MPC) as needed. Its goal is to achieve stabilization with horizons as small as possible so that MPC methods can be used on faster and/or more complicated dynamic processes. 

   \end{abstract}

\section{Introduction}
\setcounter{equation}{0}
Model Predictive Control (MPC) is a way to  steer a discrete time control system to a desired 
operating point.     We will present an extension of MPC that we call Adaptive Horizon Model Predictive Control (AHMPC) which adjusts the length of the horizon in MPC while confirming in real time that stabilization is occuring for a nonlinear system.

We are not the first to consider adaptively changing  the horizon length in MPC, see
\cite{MM93}, \cite{PY93a}, \cite{PY93b}, \cite{PY93c}.  In these papers the horizon is changed so that a terminal constraint is satisfied by the predicted state at the end of horizon.
In \cite{Gi10} the horizon length is adaptively changed to ensure that the infinite horizon cost of using the finite horizon MPC scheme is not much more than the cost of the corresponding  infinite horizon optimal control problem.

Adaptive horizon tracking is discussed in \cite{PDHW06}
and \cite{DE11}.
In \cite{KS08}  an adaptive parameter estimation algorithm  suitable for MPC was proposed, which uses the available input and output signals to estimate the unknown system parameters.
In  \cite {GPSW10} is  a
detailed analysis of the impact of the optimization horizon and the  time varying control
horizon on stability and performance of the closed loop is given.
  
  \section{Review of Model Predictive Control}
  \setcounter{equation}{0}
  We briefly describe MPC following the definitive treatise of  \cite{RM09}.
  We largely follow their notation.

  We are given a controlled, nonlinear dynamics in discrete time 
  \bea \label{dyn}
  x^+&=& f(x,u)
  \eea
  where the state $x\in\Real^{n\times1}$, the control $u\in \Real^{m\times1}$ and $x^+(k)=x(k+1)$.
 Typically this a   discretization of a controlled, nonlinear dynamics in continuous time.  The goal is 
  to find a feedback law $u(k)=\kappa(x(k))$ that drives the state of the system to some desired 
  operating point.  A pair $(x^e,u^e)$ is  an operating point  if $f(x^e,u^e)=x^e$.   We conveniently  assume that, after  state and control coordinate translations, the operating point of interest is $(x^e,u^e)=(0,0)$.    
  
  The controlled dynamics may be subject to constraints such
  as
  \bea \label{sc}
  x&\in& \mathbb{X}\subset \Real^{n\times 1}\\
  u&\in& \mathbb{U}\subset \Real^{m\times 1} \label{scu}
  \eea
  and possibly constraints involving both the state and control
  \bea \label{scc}
  y=h(x,u)&\in& \mathbb{Y}\subset \Real^{p\times 1}
  \eea
  A control $u$ is said to be feasible at $x\in \mathbb{X}$ if  
   \bea \label{feas}
   u\in \mathbb{U},&
  f(x,u) \in\mathbb{X},&
   h(x,u)\in\mathbb{Y}\label{feasxu}
   \eea
  Of course  the stabilizing  feedback $\kappa(x)$ that we seek needs to be feasible, that is, for every  $x\in\mathbb{X}$,
  \bean
  \kappa(x)\in  \mathbb{U},&
  f(x,\kappa(x))  \in \mathbb{X},&
   h(x,\kappa(x))  \in \mathbb{Y}
   \eean

  An ideal way to find a stabilizing feedback is to choose a Lagrangian $l(x,u)$ (aka running cost)  that is nonnegative definite 
  in $x,u$, positive definite in $u$ and then to solve the infinite horizon optimal control problem of minimzing   the quantity
$
 \sum_{k=0}^\infty l(x(k),u(k))
$  over  all choices of infinite control sequences \\${\bf u} =(u(0), u(1),\ldots)$  subject to the dynamics (\ref{dyn}), the constraints (\ref{sc}, \ref{scu}, \ref{scc}) and the initial condition $x(0)=x^0$.   Assuming the minimum exists for each $x^0\in\mathbb{X}$, we define the optimal cost function
    \bea \label{icost}
 V (x^0)=\min_{\bf u} \sum_{k=0}^\infty l(x(k),u(k))
  \eea
 Let ${\bf u}^* =(u^* (0), u^* (1),\ldots)$ be a minimzing control sequence with corresponding state sequence ${\bf x}^* =(x^* (0)=x^0, x^* (1),\ldots)$.  Minimizing control and state sequences need not be unique but we shall generally ignore this   problem because we are using optimization as a path to stabilization.  The key question is whether the possibly nonunique solution is stabilizing to the desired operating point.  As we shall see AHMPC nearly verifies stabilization in real time.
  
  If a pair $V (x)\in \Real, \kappa (x)\in \Real^{m\times 1}$ of functions satisfy the infinite horizon Bellman Dynamic Program Equations (BDP)
  \eq \label{BDP}
  V (x)&=&\mbox{min}_u  \left\{V (f(x,u))+ l(x,u)\right\}\\
  \kappa (x)&=& \mbox{argmin}_u \left\{ V (f(x,u))+l(x,u)\right\}   \\
\ee
and the constraints 
\bea
\kappa (x)\in \mathbb{U},&
f(x,\kappa (x))\in\mathbb{X},&
h(x,\kappa (x))\in \mathbb{Y}
 \label{cc}
\eea
for all $x\in \mathbb{X}$
then it is not hard to show that $V (x)$ is the optimal cost and
$ \kappa (x)$ is an optimal feedback 
 on $ \mathbb{X}$.
Under suitable conditions  a Lyapunov argument can be used to show that the feedback
$ \kappa (x)$ is stabilizing.

The difficulty with this approach is that it is generally impossible to  solve the BDP equations
on a large domain $\mathbb{X}$ if the state dimension $n$ is greater than $2$ or $3$.  So both theorists and practicioners have turned to Model Predictive Control (MPC). 

In MPC one chooses a Lagrangian $l(x,u)$, a horizon length $N$, a terminal domain 
$\mathbb{X}_f\subset\mathbb{X}$ containing $x=0$ and a terminal cost $V_f(x)$ defined and positive definite on  $\mathbb{X}_f$.  Then one considers the problem of minimizing 
\bean
\sum_{k=0}^{N-1} l(x(k),u(k)) + V_f(x(N))
\eean
 by choice of feasible
 \bean{\bf u}_N=(u_N(0),u_N(1),\ldots,u_N(N-1))
 \eean
subject to the dynamics (\ref{dyn}), the constraints (\ref{cc}), the final condition $x(N)\in \mathbb{X}_f$ and the initial condition $x(0)=x^0$.
Assuming this problem is solvable, let $V_N(x^0)$ denote the optimal cost, 
\bea  \label{Ncost}
V_N(x^0)&=&\min_{{\bf u}_N} \sum_{k=0}^{N-1} l(x(k),u(k)) + V_f(x(N))
\eea
and let
\bean
{\bf u}_N^*&=&(u_N^*(0),u_N^*(1),\ldots,u_N^*(N-1))\\ {\bf x}_N^*&=&(x_N^*(0)=x^0,x_N^*(1),\ldots,x_N^*(N))
\eean denote optimal control and state sequences when the horizon length is $N$.  We then  define the MPC feedback law
$
\kappa_N(x^0)= u_N^*(0)
$.

The terminal set $\mathbb{X}_f$ is controlled  invariant (aka viable) if for each $x\in\mathbb{X}_f$
there exists a $u\in \mathbb{U}$ such that $f(x,u)\in \mathbb{X}_f$ and the constraint  (\ref{scc})
is satisfied.   If $V_f(x)$ is a control Lyapunov function on $\mathbb{X}_f$ then, under suitable conditions, a Lyapunov argument can be used to show that the feedback
$ \kappa_N(x)$ is stabilizing on $\mathbb{X}_f$.   See \cite{RM09} for more details.

AHMPC requires a little more, the existence of terminal feedback $u=\kappa_f(x)$ defined on the terminal set $\mathbb{X}_f$ that leaves it positively invariant, if $x\in \mathbb{X}_f$ then 
$f(x,\kappa_f(x)) \in \mathbb{X}_f$, and
 which makes $V_f(x)$ a strict  Lyapunov function on $\mathbb{X}_f$ for the closed loop dynamics, if $x \in \mathbb{X}_f$ and $x\ne 0$ then
\bean
V_f(x)>V_f(f(x,\kappa_f(x)))\ge 0 
\eean

If $x\in \mathbb{X}_f$ then AHMPC does not need to solve (\ref{Ncost}) to get $u$, it just takes $u=\kappa_f(x)$.   A similar scheme has been called dual mode control in \cite{MM93}.

The advantage of solving the  finite  horizon   optimal control problem (\ref{Ncost}) over solving the infinite horizon problem (\ref{icost}) is that
it may be possible to solve the former  on-line as the process evolves.  If it is known that the terminal set $\mathbb{X}_f$ can be reached from 
the current  state  $x$ in $N$ or fewer steps then  the finite horizon $N$ optimal control problem is a feasible nonlinear program with finite dimenionsal decision variable ${\bf u}_N\in \Real^{m\times N}$.   If the time step is long enough, if $f,h,l$ are reasonably simple and if $N$ is small enough then this nonlinear program possibly can be solved in a fraction of one time step for ${\bf u}^*_N$.  Then the first element of this sequence $u_N^*(0)$ is used as the control at the current time.  The system evolves one time step and the process is repeated at the next time.  Conceptually MPC computes a feedback law $\kappa_N(x) =u^*_N(0)$  but only  at  values of $x$ when and where it is needed.

Some authors like \cite{GMTT05},  \cite{Gr12} do away with the terminal cost $V_f(x)$  but there is a theoretical reason and  a practical reason to use one.  The theoretical reason is that a control Lyapunov terminal cost facilitates a proof of asymptotic stability via a simple Lyapunov argument, see \cite{RM09}. But this is not  a binding reason because under suitable assumptions,  asymptotic stability can be shown even when there is no terminal cost provided the horizon is sufficiently long.  The practical reason is more important, when there is a terminal cost one can usually use a shorter horizon $N$.  A shorter horizon reduces the dimension $mN$ of the decision variables in the nonlinear programs  that need to be solved on-line.  Therefore MPC with a suitable terminal cost can be used for faster and more complicated systems. 

 An ideal terminal cost $V_f(x)$ is $V (x)$ of the corresponding infinite horizon optimal control  provided that the latter can be accurately computed off-line  on a reasonably large  terminal set $\mathbb{X}_f$.   For then the  infinite horizon cost (\ref{icost}) and (\ref{Ncost}) will be the same.   One should not make too much of this fact as stabilization is our goal, the optimal control problems are just a means to accomplish this.   This in contrast to Economic MPC where the cost and the associated  Lagrangian are  chosen to model real world costs.

\section{ Adaptive Horizon Model Predictive Control}
For AHMPC we  assume that we have the following.   
\begin{itemize}
 \item A discrete time dynamics $f(x,u)$ with operating point $x=0,u=0$.
 \item A Lagrangian  $l(x,u)$, nonegative definite in $(x,u)$ and positive definite in $u$.
\item State  constraints $x\in \mathbb{X}$ where $\mathbb{X}$ is a  neighborhood of $x=0$.
\item Control  constraints $u\in \mathbb{U}$ where $\mathbb{U}$ is a  neighborhood of $u=0$.
\item Mixed constraints $ h(x,u)\in  \mathbb{Y}$ which are not active at the operating point $x=0,u=0$.
\item The dynamics is recursively feasible on $ \mathbb{X}$, that is, for every $x\in \mathbb{X}$ there
is a $u\in \mathbb{U}$ satisfying $ h(x,u)\in  \mathbb{Y}$ and $f(x,u)\in  \mathbb{X}$
\item A terminal cost $V_f(x)$ defined and nonnegative definite on  some neighborhood  $\mathbb{X}_f$  of the operating point $x=0,u=0$.
The neighborhood $\mathbb{X}_f$
need not be known explicitly.
\item A terminal feedback $u=\kappa_f(x)$ defined on  $\mathbb{X}_f$ such that
the terminal cost is a valid Lyapunov function on $\mathbb{X}_f$ for the closed loop dynamics
using the terminal feedback $u=\kappa_f(x)$.
\end{itemize}

One way of obtaining a terminal pair $V_f(x), \ \kappa_f(x)$ is to approximately solve the infinite horizon dynamic program equations (BDP) on some neighborhood of the origin.  For example if the linear part of the dynamics and the quadratic part of the Lagrangian constitute a  LQR problem satisfying the standard conditions then  one can let 
$V_f(x)$ be the quadratic optimal cost and $\kappa_f(x)$ be the linear optimal feedback of the LQR problem.    Of course the problem with such terminal pairs $V_f(x), \ \kappa_f(x)$ is that generally there is no way to estimate the terminal set $\mathbb{X}_f$ on which the feasibility and Lyapunov conditions are satisfied.  It is reasonable to expect that they are satisfied on some terminal set but the extent of this terminal set is very difficult to estimate.

In the next section  we show how higher degree Taylor polynomials for the optimal cost and optimal feedback can be computed by the discrete time version of  \cite{Al61} found in \cite{AK12} because this can lead to a larger  terminal set $\mathbb{X}_f$ on which the feasibility and the Lyapunov conditions are satisfied.    It would be very difficult to determine what this terminal set is but we do not need to do this.

AHMPC mitigates this last difficulty just as MPC mitigates the problem of solving the infinite horizon Bellman Dynamic Programming equations BDP (\ref{BDP}).  MPC does not try to compute the optimal cost and optimal feedback everywhere, instead it computes them just when and where they are needed.  AHMPC does not try to compute  the set $\mathbb{X}_f$ on which  $\kappa_f(x)$ is feasible and stabilizing, it just tries to determine if the  end state $x_N^*(N)$ of the currently computed optimal trajectory  is in a terminal  set $\mathbb{X}_f$ where  the feasibility and Lyapunov conditions are satisfied.  
 
 Suppose the current state is $x$ and we have solved the horizon $N$ optimal control problem for
$ {\bf u}^*_N=(u^*_N(0),\ldots, u^*_N(N-1))$,
${\bf x}^*_N=(x^*_N(0)=x,\ldots, x^*_N(N))$
  The terminal 
 feedback $u=\kappa_f(x)$ is used to compute  $M$ additional steps of this state trajectory
  \bea \label{M}
 x^*_N( k+1))&=& f(x^*_N(k),\kappa_f(x^*_N(k))
 \eea 
for $k=N,\ldots, N+M-1$. 

Then one checks that the feasibility and Lyapunov conditions hold for the extended part of the state sequence,
\bea
\label{cf}
\kappa_f(x^*_N(k))&\in & \mathbb{U}\\
f(x^*_N(k),\kappa_f(x^*_N(k)))&\in & \mathbb{X}  \label{sf}\\
h(x^*_N(k),\kappa_f(x^*_N(k)))&\in & \mathbb{Y}  \label{scf}\\
 \label{L1}
V_f (x^*_N(k))&\ge & \alpha(|x^*_N(k)|)\\
V_f (x^*_N(k))-V_f (x^*_N(k+1))&\ge & \alpha(|x^*_N(k)|)
\label{L2}
\eea
for $k=N,\ldots, N+M-1$ and some Class K function $\alpha(s)$.  For more on Class K functions we refer the reader to \cite{Kh96}. 

If (\ref{cf}-\ref{L2}) hold for all for $k=N,\ldots, N+M-1$.  then we presume that $x^*_N(N)\in \mathbb{X}_f$, a set  where $\kappa_f(x)$ is stabilizing
and we use the control $u_N^*(0)$ to move one time step forward to  $x^+=f(x, u^*_N(0))$.  At this next state $x^+$ we 
solve the horizon $N-1$ optimal control problem and check that the extension of the new optimal trajectory satisfies (\ref{cf}-\ref{L2}).

If (\ref{cf}-\ref{L2}) does not hold for all for $k=N,\ldots, N+M-1$.  then we presume that $x^*_N(N)\notin \mathbb{X}_f$.
We extend current horizon by $L\ge1$  and if time permits we solve the horizon $N+L$ optimal control problem at the current state $x$ and then check the feasibility and Lyapunov  conditions again.  We keep increasing $N$ by $L$   until these conditions are satisfied on the extension of the trajectory.
If we run out of time before  they are satisfied then we use the last computed $u_N^*(0)$ and move one time step forward to  $x^+=f(x, u^*_N(0))$.  At $x^+$  we 
solve the horizon $N+L$ optimal control problem.

How does one choose the extended horizon $M$ and the class ${\it K}$ function $\alpha(\cdot)$?   If the extended part of the state sequence is actually in the region where the terminal cost $V_f (x)$ and the terminal feedback $\kappa_f(x)$ well approximate the solution to infinite horizon optimal control problem then the   dynamic programing equations (\ref{BDP}) should approximately hold.  In other words  
\bean
V_f (x^*_N(k))-V_f (x^*_N(k+1))&\approx& l(x^*_N(k), \kappa_f(x^*_N(k))\ \ge \ 0
\eean
If this does not hold throughout the extended trajectory we should increase both  the horizon $N$. We can also increase the extended horizon $M$ but this may not be necessary.   If the Lyapunov and feasibilty conditions are going to fail somewhere on the extension
it is most likely this will happen at the beginning of the extension.
Also we should choose $\alpha(\cdot)$ so that $\alpha(|x|)<|l(x,\kappa_f(x))|/2$.  

The nonlinear programming problems generated by employing MPC on a nonlinear system are generally nonconvex so the solver might return local rather than global minimizers.   In which case there is no guarantee that an MPC approach is actually stabilizing.   AHMPC mitigates this difficulty by checking that stabilization is occuring.   If  (\ref{cf}-\ref{L2}) don't hold even after the horizon $N$ has been increased substantially then this is a strong indication that the solver is returning locally rather than globally minimizing  solutions and these local solutions are not stabilizing.   To change this behaviour one needs to start the solver at a substantially different initial guess.   Just how one does this is an open research question.   It is essentially the same question as which initial guess should one pass to the solver at the first step of MPC.

  The actual computation of the $M$ additional steps
(\ref{M}) can be done very quickly because the closed loop dynamics function $f(x,\kappa_f(x))$ can be computed and compiled before hand.  Similarly the feasibility and Lyapunov conditions (\ref{cf}-\ref{L2}) can be computed and compiled before hand. The number $M$ of additional time steps is a design parameter.  One choice is to  take $M$ a fraction of the current $N$.

\section{Choosing a Terminal Cost and a Terminal Feedback }
 A standard way of obtaining a terminal cost $V_f (x)$ and a terminal feedback $\kappa_f(x)$ is to solve the Linear Quadratic Regulator (LQR) using the quadratic part of the Lagrangian and the linear part of dynamics around the operating point $(x^e,u^e)=(0,0)$.    
Suppose
\bean
l(x,u)&=& {1\over 2} \left(x'Qx+2x'Su+u'Ru\right) +O(x,u)^3\\
f(x,u)&=&Fx+Gu+O(x,u)^2
\eean
Then the LQR problem is to find $P,K$  such that 
 \bea 
P &=&F'P  F-\left(F'P  G+S\right)\left(R+G'P  G\right)^{-1}\left(G'P  F+S'\right)+Q\nonumber\\
\label{P}\\
K &=& -\left(R+G'P  G\right)^{-1}\left(G'P  F+S'\right) \label{K}
\eea
 Under mild assumptions, the stabilizability of $F,G$, the detectability of $Q^{1/2},F$, the nonnegative definiteness of $[Q,S;S',R]$ and the positive definiteness of $R$,  there exist an unique nonnegative definite $P $ satisfying the first equation (\ref{P}) which is called the discrete time algebraic Riccati equation.  Then  $K $ given by (\ref{K}) puts all the poles of the closed loop linear dynamics
 \bean
 x^+&\left(F+GK\right)x
 \eean 
 inside of the open unit disk.  See \cite{AM97} for details.  But $l(x,u)$ is a design parameter so if we  choose $[Q,S;S',R]$ to be positive definite and then $P$ will be positive definite.

 If we define the terminal cost to be $V_f (x)={1\over 2} x'P x $ then we know that it is positive definite for all $x\in \Real^{n\times 1}$.   If we define the terminal feedback to be $\kappa_f(x)=Kx$ then we know by a Lyapunov argument that the nonlinear closed loop dynamics
 \bean
 x^+&=&f(x,\kappa(x))
 \eean
  is locally asymptotically stable around $x^e=0$.  The problem is that we don't know what is the neighborhood $\mathbb{X}_f$ of asymptotic stability and computing it off-line can be very difficult in state dimensions higher that two or three.
  
There are other possible choices for the terminal cost and terminal feedback.  Al'brekht \cite{Al61} showed how the Taylor polynomials of the optimal cost and the optimal feedback could be computed for some smooth, infinite horizon optimal control problems in continuous time.   Aguilar and Krener \cite{AK12} extended this to some smooth, infinite horizon optimal control problems in discrete  time.  The discrete time Taylor polynomials of the optimal cost and the optimal feedback may be used as the terminal cost and terminal feedback in an AHMPC scheme so we briefly review \cite{AK12}.

Since we assumed that $f(x,u)$ and $l(x,u)$ are smooth and the constraints are not active at the origin, we can simplify the BDP equations.
The simplified Bellman Dynamic Programming equations (sBDP) are obtained by setting the derivative with respect to $u$ of the quantity to be minimized in (\ref{BDP}) to zero.  The result is 
 \bea \label{sBDP1}
  V (x)&=&V (f(x,\kappa (x))+ l(x,\kappa (x))\\ \nonumber
  \\
  0&=& \frac{\partial V }{ \partial x}(f(x,u))\frac{\partial f}{\partial u}(x,  \kappa (x) ) +\frac{\partial l}{\partial u}(x,\kappa (x)) 
\label{sBDP2}
\eea
If  the quantity to be minimized is strictly convex in $u$ then the BDP equations and sBDP equations are equivalent.  But if not then BDP implies sBDP but not vice versa.

Suppose the discrete time dynamics and Lagrangian have Taylor polynomials around the operating point $x=0,u=0$ of the form
\bean
f(x,u)&=& Fx+Gu+f^{[2]}(x,u)+\cdots+f^{[d]}(x,u)\\
l(x,u)&=&{1\over 2} \left(x'Qx+2x'Su+u'Ru\right)+l^{[3]}(x,u)+\cdots+ l^{[d+1]}(x,u)
 \eean
for some integer $d\ge 1$ and where $^{[j]}$ indicates homogeneous polynomial terms of degree $j$.

Also suppose the infinite horizon optimal cost and optimal feedback have similar Taylor polynomials 
\bean
V (x)&=&{1\over 2} x'P  x+V ^{[3]}(x)+\cdots+ V ^{[d+1]}(x)\\
\kappa (x)&=& K  x+\kappa ^{[2]}(x,u)+\cdots+\kappa ^{[d]}(x,u)
 \eean

We plug these polynomials into sBDP and collect terms of lowest degree.  The lowest degree in (\ref{sBDP1}) is two while 
in (\ref{sBDP2}) it is one.  The result is the discrete time Ricatti equations (\ref{P}, \ref{K}).

At the next degrees we obtain the equations
\bea
V^{[3]}(x) - V^{[3]}((F+GK)x) &=& ((F+GK)x)'Pf^{[2]}(x,Kx)\nonumber\\&&+l^{[3]}(x,Kx)\label{BDPa}\\
\left(\kappa^{[2]}(x)\right)'\left(R+G'PG\right)&=& -\frac{\partial  V^{[3]}}{\partial x}((F+GK)x)G\nonumber\\
&&-((F+GK)x)'P\frac{\partial f^{[2]}}{\partial u}(x,Kx)\nonumber\\&&
-\frac{\partial l^{[3]}}{\partial u}(x,Kx) \label{BDPb}
\eea

Notice these are linear equations in the unknowns $V^{[3]}(x)$ and $\kappa^{[2]}(x)$ and the right sides of these equations involve only known quantities.   Moreover  $\kappa^{[2]}(x)$
does not appear in the first equaion.
The eigenvalues of the linear operator 
\bea \label{op3}
V^{[3]}(x) &\mapsto&V^{[3]}(x) - V^{[3]}((F+GK)x) 
\eea
are of the form of $1$ minus the products of three eigenvalues of $F+GK$.  Since the 
 eigenvalues of $F+GK$ are all inside the open unit disk this operator (\ref{op3})  is invertible.
 Having solved (\ref{BDPa}) for $V^{[3]}(x)$ we can readily solve (\ref{BDPb}) 
 for $\kappa^{[2]}(x)$ since we have assumed $R$ is positive definite
 
 The higher degree equations are similar, at degrees $d+1, d$ they take the form 
 \bean
 V^{[j+1]}(x) - V^{[j+1]}((F+GK)x) &=& \mbox{ Known Quantities}\\
\left(\kappa^{[j]}(x)\right)'\left(R+G'PG\right)&=&\mbox{ Known Quantities} 
\eean 
 The "Known Quantities" involve the terms of the Taylor polynomials of $f,l$ and previously computed
 $V^{[i+1]}, \kappa^{[i]}$ for $1\le i<j$.
 Again the equations are linear in the unknowns  $V^{[k+1]}, \kappa^{[k]}$ and the first equation
 does not involve $\kappa^{[k]}$.   The eigenvalues of the linear operator 
 \bean
V^{[k+1]}(x) &\mapsto&V^{[k+1]}(x) - V^{[k+1]}((F+GK)x) 
\eean
 are of the form of $1$ minus the products of $k+1$ eigenvalues of $F+GK$. 
 We have  written MATLAB code to solve these equations to any degree and in any dimensions.  The code is quite fast.
It found the Taylor polynomials of the optimal cost to degree 6 and the optimal feedback to degree 5 in 0.12 sec. on a laptop
 using 3.1 GHz Intel Core i5.

\section{Completing the Squares} 
 The infinite horizon optimal cost is certainly nonnegative definite and if we choose $Q>0$
 then it is positive definite.   That implies that its quadratic part $V^{[2]}(x)={1\over 2}x'Px$
 is positive definite.
 
 But its Taylor polynomial of degree $d+1$,
 \bean
 V^{[2]}(x)+V^{[3]}(x)+\cdots+V^{[d+1]}(x)
 \eean 
 need not be positive definite for $d>1$.  This can lead to problems
 if we define this Taylor polynomial to be our terminal cost $V_f(x)$ because then the nonlinear program solver
 might return a negative cost $V_N(x)$ (\ref{Ncost}).
 The way around this difficulty is to "complete the squares".
 
 {\bf Theorem}
Suppose  a polynomial $V(x)$ is of degrees two through $d+1$ in $n$ variables $x_1,\ldots,x_n$. 
If the quadratic part of $V(x)$ is positive definite then there exist a nonnegative definite polynomial $W(x)$ of degrees two through
$2d$ such that the part of $W(x)$ that is of degrees two through $d+1$ equals $V(x)$.   Moreover, we know that $W(x)$ is nonegative definite because it is the sum of $n$ squares.
\\

{\bf Proof:}
We start with the quadratic part of $V(x)$, because it is positive definite it must be of the form ${1\over 2}x'Px$ where $P$ is a positive definite $n\times n$ matrix.  We know that there is an orthogonal matrix $T$ that diagonalizes $P$
\bean
T'PT=\bmt \lambda_1& & 0\\&\ddots &\\0&& \lambda_n \emt
\eean
where $\lambda_1\ge \lambda _2\ge \ldots \ge \lambda_n>0$.  We make the linear change of coordinates $x=Tz$.  We shall show that
$V(z)=V(Tz)$ can be extended with higher degree terms to a polynomial  $W(z)$ of degrees two through $2d$ which is a sum of $n$ squares.
  We do this degree by degree.  We have already showed that the degree two part of $V(z)$ is a sum of $n$ squares,
  \bean
{1\over 2} \sum_{i+1}^n \lambda_i z_i^2&=&{1\over 2}\left( \lambda_1 z_1^2+\cdots  \lambda_n z_n^2\right)
  \eean
  The degrees two and three parts of $V(z)$ are of the form
  \bean
{1\over 2} \sum_{i+1}^n \lambda_i z_i^2+\sum_{i_1=1}^n\sum_{i_2=i_1}^n\sum_{i_3=i_2}^n \gamma_{i_1,i_2,i_3}z_{i_1}z_{i_2}z_{i_3}
  \eean
  
  Consider the expression
  \bean
\delta_1(z)&=& z_1 +  \sum_{i_2=1}^n\sum_{i_3=i_2}^n \delta_{1,i_2,i_3}z_{i_2}z_{i_3}
\eean
then
\bean
{\lambda_1\over 2} (\delta_1(z))^2&=& {\lambda_1\over 2}\left(z_1 +  \sum_{i_2=1}^n\sum_{i_3=i_2}^n \delta_{1,i_2,i_3}z_{i_2}z_{i_3}\right)^2\\
&=& {\lambda_1 \over 2} \left( z_1^2
 +2 \sum_{i_2=1}^n\sum_{i_3=i_2}^n \delta_{1,i_2,i_3}z_1z_{i_2}z_{i_3}+{1\over 2}( \sum_{i_2=1}^n\sum_{i_3=i_2}^n \delta_{1,i_2,i_3}z_{i_2}z_{i_3} )^2\right)
  \eean

Let $\delta_{1,i_2,i_3}= {\gamma_{1,i_2.i_3}\over \lambda_1}$ then the degrees two and  three parts of 
\bean \label{pp1}
V(z)-{\lambda_1\over 2}(\delta_1(z))^2
\eean
have no terms involving $z_1$.

Next consider the expression
  \bean \label{z2}
\delta_2(z)&=& z_2 +  \sum_{i_2=2}^n\sum_{i_3=i_2}^n \delta_{2,i_2,i_3}z_{i_2}z_{i_3} 
\eean
then 
\bean
{\lambda_2\over 2} (\delta_2(z))^2
&=& {\lambda_2 \over 2} \left( z_2^2
 +2 \sum_{i_2=2}^n\sum_{i_3=i_2}^n \delta_{2,i_2,i_3}z_2z_{i_2}z_{i_3}+{1\over 2}( \sum_{i_2=2}^n\sum_{i_3=i_2}^n \delta_{2,i_2,i_3}z_{i_2}z_{i_3} )^2\right)
  \eean
Let $\delta_{2,i_2,i_3}= {\gamma_{2,i_2,i_3}\over \lambda_2}$ then the degrees two and  three parts of 
\bean
V(z)-{\lambda_1\over 2}(\delta_1(z))^2-{\lambda_2\over 2}(\delta_2(z))^2
\eean
have no terms involving either $z_1$ or $z_2$.

We continue on in this fashion defining $\delta_3(z), \ldots, \delta_n(z)$ such that 
\bean
V(z)-\sum_{i=1}^n{\lambda_i\over 2}(\delta_i(z))^2&=&\sum_{i_1=1}^n\sum_{i_2=i_1}^n\sum_{i_3=i_2}^n \sum_{i_4=i_3}^n  \gamma_{2,i_2,i_3,i_4}z_{i_1}z_{i_2}z_{i_3}z_{i_4}+O(z)^5
\eean
has no terms of  degrees either two or  three.

We redefine
 \bean
\delta_1(z)&=& z_1 +  \sum_{i_2=1}^n\sum_{i_3=i_2}^n \delta_{1,i_2,i_3}z_{i_2}z_{i_3}+  \sum_{i_2=1}^n\sum_{i_3=i_2}^n\sum_{i_4=i_3}^n \delta_{1,i_2,i_3,i_4}z_{i_2}z_{i_3}z_{i_4} 
\eean 

This does not change the degree two and three terms of ${\lambda_1\over 2}(\delta_1(z))^2$ and its degree four terms are of the form
\bean
\lambda_1\ \sum_{i_2=1}^n\sum_{i_3=i_2}^n\sum_{i_4=i_3}^n \delta_{1,i_2,i_3,i_4}z_1z_{i_2}z_{i_3}z_{i_4} 
\eean
If we let $\delta_{1,i_2,i_3,i_4}={ \gamma_{1,i_2,i_3,i_4}\over \lambda_1}$  then we cancel the degree four terms involving $z_1$ in
\bean
V(z)-\sum_{j=1}^n{\lambda_j\over 2}(\delta_j(z))^2
\eean 

Next we redefine 
\bean
\delta_2(z)&=&z_2 +  \sum_{i_2=2}^n\sum_{i_3=i_2}^n \delta_{2,i_2,i_3}z_{i_2}z_{i_3}+  \sum_{i_2=2}^n\sum_{i_3=i_2}^n\sum_{i_4=i_3}^n \delta_{2,i_2,i_3,i_4}z_{i_2}z_{i_3}z_{i_4} 
\eean 
Again this does not change the degree two and three terms of ${\lambda_2\over 2}(\delta_2(z))^2$ and its degree four terms are of the form
\bean
\lambda_2\ \sum_{i_2=2}^n\sum_{i_3=i_2}^n\sum_{i_4=i_3}^n \delta_{2,i_2,i_3,i_4}z_2z_{i_2}z_{i_3}z_{i_4} 
\eean
If we let $\delta_{2,i_2,i_3,i_4}={ \gamma_{2,i_2,i_3,i_4}\over \lambda_2}$  then we cancel the degree four terms involving  $z_2$ in
\bean
V(z)-\sum_{j=1}^n{\lambda_j \over2}(\delta_j(z))^2
\eean 

We continue on in this fashion.  The result is a  sum of squares whose degree two through four terms equal $V(z)$.

Eventually we define
\bean
\delta_j(z)&=& z_j+  \sum_{i_2=j}^n\sum_{i_3=i_2}^n \delta_{j,i_2,i_3}z_{i_2}z_{i_3}+\cdots\\
&&+  \sum_{i_2=j}^n \ldots \sum_{i_d=i_{d-1}}^n\delta_{j,i_2,\ldots,i_d}z_{i_2}\cdots z_{i_d} 
\eean
and
\bean
W(z)&=& \sum_{j=1}^n {\lambda_j\over 2}\left( \delta_j(z)\right)^2
\eean
  
At degree $d=3$ we solved a linear equation from the quadratic coefficients of $\delta_1(z),\ldots, \delta_n(z)$ to the cubic coefficients of $V(z)$.  We restricted the domain of this mapping by requiring that the quadratic part of $\delta_i(z)$ does not depend on $z_1,\ldots,z_{i-1}$.  This made the restricted mapping square, the dimensions of the domain  and the range  of the linear mapping are the same.  We showed that the restricted mapping has an  unique solution.  If we drop this restriction then the overall dimension of the domain is
\bean
n\left(\begin{array}{ccc} n+d-1\\ d\end{array}\right)
\eean
while the dimension of the range is
\bean
\left(\begin{array}{ccc} n+d\\ d+1\end{array}\right)
\eean
So the  unrestricted mapping has more unkowns than equations and hence there are multiple solutions.

But the restricted solution that we constructed is a least squares solution to the unrestricted equations because $\lambda_1\ge \lambda _2\ge \ldots \ge \lambda_n>0$.  To see this consider the coefficient $\gamma_{1,1,2}$ of $z_1^2z_2$.  If we allow $\gamma_2(z)$ to have a term of the form $\delta_{2,1,1} z_1^2$ we can also cancel $\gamma_{1,1,2}$ by choosing $\delta_{1,1,2}$ and $\delta_{2,1,1}$ so that 
\bean
\gamma_{1,1,2}&=&  \lambda_1\delta_{1,1,2}+ \lambda_2\delta_{2,1,1}
\eean
Because $\lambda_1\ge \lambda _2$ a least squares solution to this equation is $\delta_{1,1,2}={\gamma_{1,1,2}\over  \lambda_1}$ and $\delta_{2,1,1}=0$.
Because $T$ is orthogonal $W(x)=W(T'z)$ is also a least squares solution.

\section{Example}
Consider a double pendulum that we wish to stabilize to straight up.   The first two states are the the angles between the two links and straight up  measured in radians counterclockwise.  The other  two states are their  angular velocities.   The controls are the torques applied at the base of the lower link and at the joint between the links.  The links are assumed to massless, the base link is one meter long and the other link is two meters long.  There is a mass of two kilograms at the joint between the links and a mass of one kilogram at the tip of the upper link.  There is linear damping at the base and the joint with both coefficients equal to $0.5$ sec$^{-1}$.   The resulting  continuous time dynamics is discretized using Euler's method with time step $0.1$ sec.

We simulated AHMPC with two different terminal costs and terminal feedbacks.  
 In both cases the Lagrangian was 
\bea  \label{Lag}
{0.1\over 2}\left(|x|^2+|u|^2\right)
\eea

The first pair $V^2_f(x), \kappa^1_f(x)$  was found by solving the infinite horizon LQR problem obtained by taking the linear part of the dynamics around the operating point $x=0$ and the quadratic Lagrangian (\ref{Lag}).   Then $V^2_f(x)$ is a positive definite quadratic and $\kappa^1_f(x)$ is linear.

The second pair $V^6_f(x), \kappa^5_f(x)$  was found using the discrete time version of Al'brekht's method.  Then $V^6_f(x)$ is
the Taylor polynomial of the optimal cost to degree $6$ and $\kappa^5_f(x)$ is the Taylor polynomial of the optimal feedback to degree $5$.  But $V^6_f(x)$ is not positive definite so we completed the squares to get $V^{10}_f(x)$ which is positive definite.

In all the simulations  we imposed the control constraint $|u|_\infty\le 5$ and started at $x(0)=(0.9\pi,0.9\pi,0,0)$ with an initial horizon of $N=50$ time steps.  The extended horizon was kept constant at $M=5$.   The class K function was taken to be $\alpha(s)=s^2/10$. 

 If the Lyapunov and feasibility conditions were violated the horizon $N$ was increased by $5$ and the finite horizon nonlinear program was solved again without advancing the system.  If after three tries the Lyapunov and feasibility conditions were still not satisfied then the first value of the control sequence was used, the simulation was adavnced one time step and the horizon was increased by $5$.

If the Lyapunov and feasibility conditions were satisfied over the extended horizon then the simulation was advanced one time step and the horizon $N$ was decreased by $1$.

The simulations were first run with no noise and the results are shown in the following figures.
Both methods stabilized the links to straight up in about $t=80$ times steps ($8$ seconds).
 The degree $2d=10$ terminal cost and the degree  $d=5$ terminal  feedback  seems to do it a little more smoothly
 and with shorter maximum horizon $N=65$ versus $N=75$ for LQR ($d=1$).

\begin{figure}
 \centering
\includegraphics[width=3.5in]{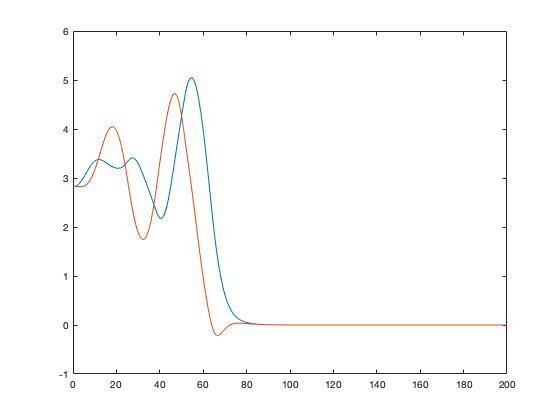}
\caption{Angles, LQR Terminal Cost and Terminal Feedback}    
\label{fig1}    
 \end{figure}
 
\begin{figure}
 \centering
\includegraphics[width=3.5in]{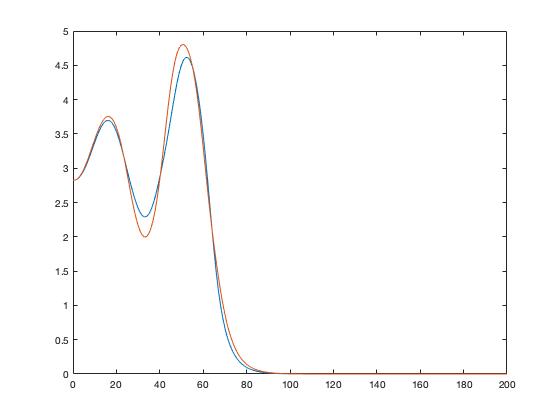}
\caption{Angles, Degree 10 Terminal Cost and Degree 5 Terminal Feedback}        
\label{fig2}
 \end{figure}
 
 \begin{figure}
 \centering
\includegraphics[width=3.5in]{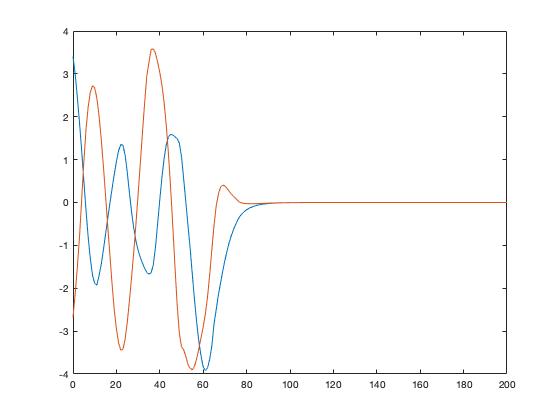}
\caption{Controls, LQR Terminal Cost and Terminal Feedback}    
\label{fig3}    
 \end{figure}
 
\begin{figure}
 \centering
\includegraphics[width=3.5in]{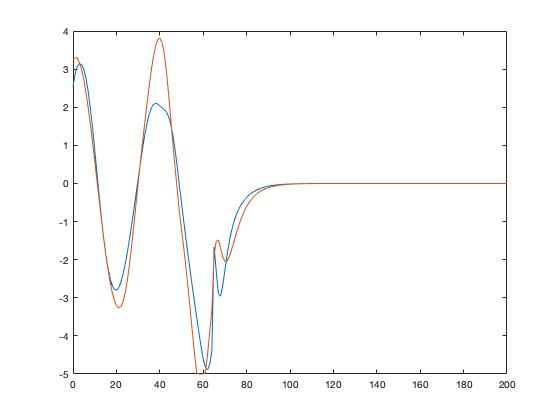}
\caption{Controls, Degree 10 Terminal Cost and Degree 5 Terminal Feedback}        
\label{fig4}
 \end{figure}

\begin{figure}
 \centering
\includegraphics[width=3.5in]{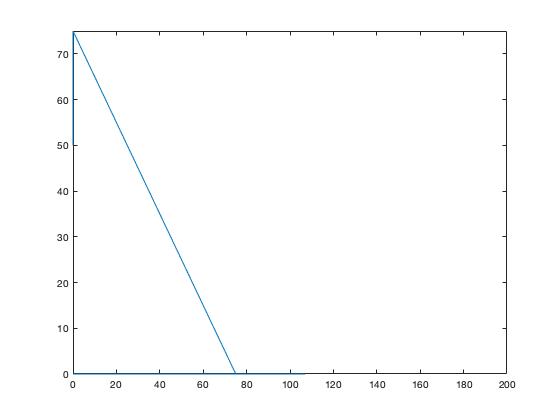}
\caption{Horizons, LQR Terminal Cost and Terminal Feedback}    
\label{fig5}    
 \end{figure}
 
\begin{figure}
 \centering
\includegraphics[width=3.5in]{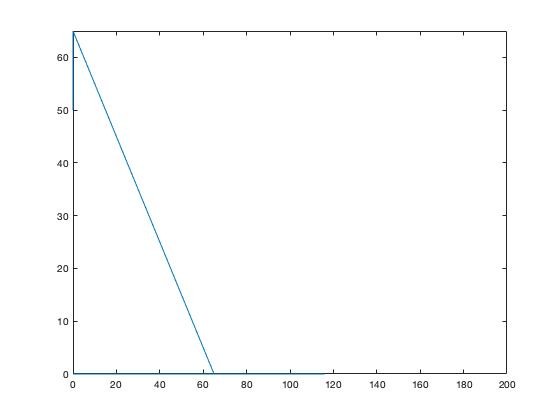}
\caption{Horizons, Degree 10 Terminal Cost and Degree 5 Terminal Feedback}        
\label{fig5}
 \end{figure}

Then we added noise to the simulations.   At each advancing step a Gaussian random vector
with mean zero and covariance $0.0004$ times the identity was added to the state.  The next figures show the results 
using the degree $10$ terminal cost and degree $5$ terminal feedback. Notice that the horizon converges to zero after $64$ time steps and stays at zero for awhile where the terminal feedback is used.  But then the noise causes the horizon to jump to $15$ before it settles back to zero.   The LQR terminal cost and feedback failed to stabilize the pendula.

\begin{figure}
 \centering
\includegraphics[width=3.5in]{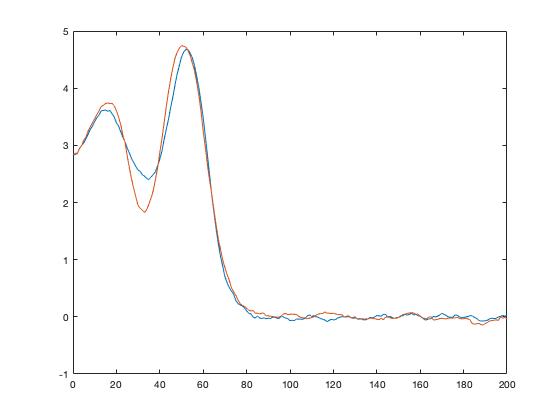}
\caption{Noisy Angles, Degree 10 Terminal Cost and Degree 5 Terminal Feedback}        
\label{fig7}
 \end{figure}
 
\begin{figure}
 \centering
\includegraphics[width=3.5in]{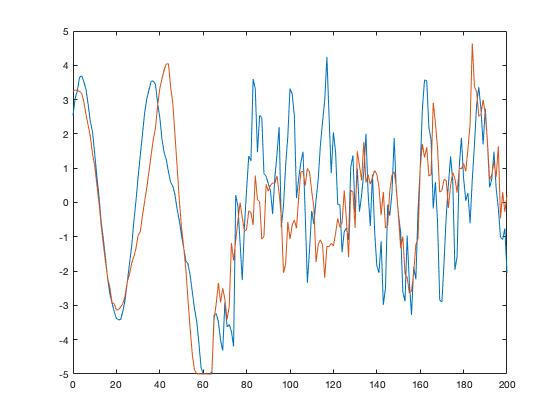}
\caption{Noisy Controls, Degree 10 Terminal Cost and Degree 5 Terminal Feedback}        
\label{fig8}
 \end{figure}
 
\begin{figure}
 \centering
\includegraphics[width=3.5in]{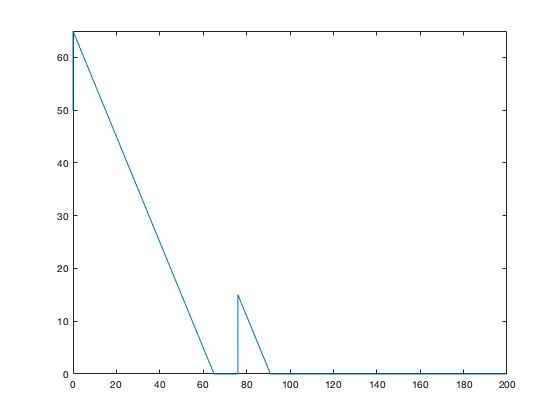}
\caption{Noisy Horizons, Degree 10 Terminal Cost and Degree 5 Terminal Feedback}        
\label{fig9}
 \end{figure}
 
We also considered $d=3$ so that after completing the squares the terminal cost is degree $2d=6$ and the terminal feedback is degree  $d=3$.  It stabilized the noiseless simulation with a maximum horizon of $N=80$ which is greater than the maximum horizons for both $d=1$ and $d=5$.  But it did not stabilize the noisy simulation.  Perhaps the reason is revealed by Taylor polynomial approximations to $\sin x$ as shown in Figure 10.   The linear approximation  in green overestimates the magnitude of $\sin x$ so the linear feedback is stronger than it needs to be to overcome gravity.
The cubic approximation  in blue underestimates  the magnitude of $\sin x$ so the cubic feedback is weaker than it needs to be to overcome gravity.  The quintic approximation in orange overestimates the magnitude of $\sin x$ so the quintic feedback is also stronger than it needs to be but by a lesser margin than the linear feedback to overcome gravity.  This may explain why the degree $5$ feedback stabilizes the noise free pendula in a smoother fashion than the linear feedback, compare Figures 1 and 2 and Figures 3 and 4.

\begin{figure}
 \centering
\includegraphics[width=3.5in]{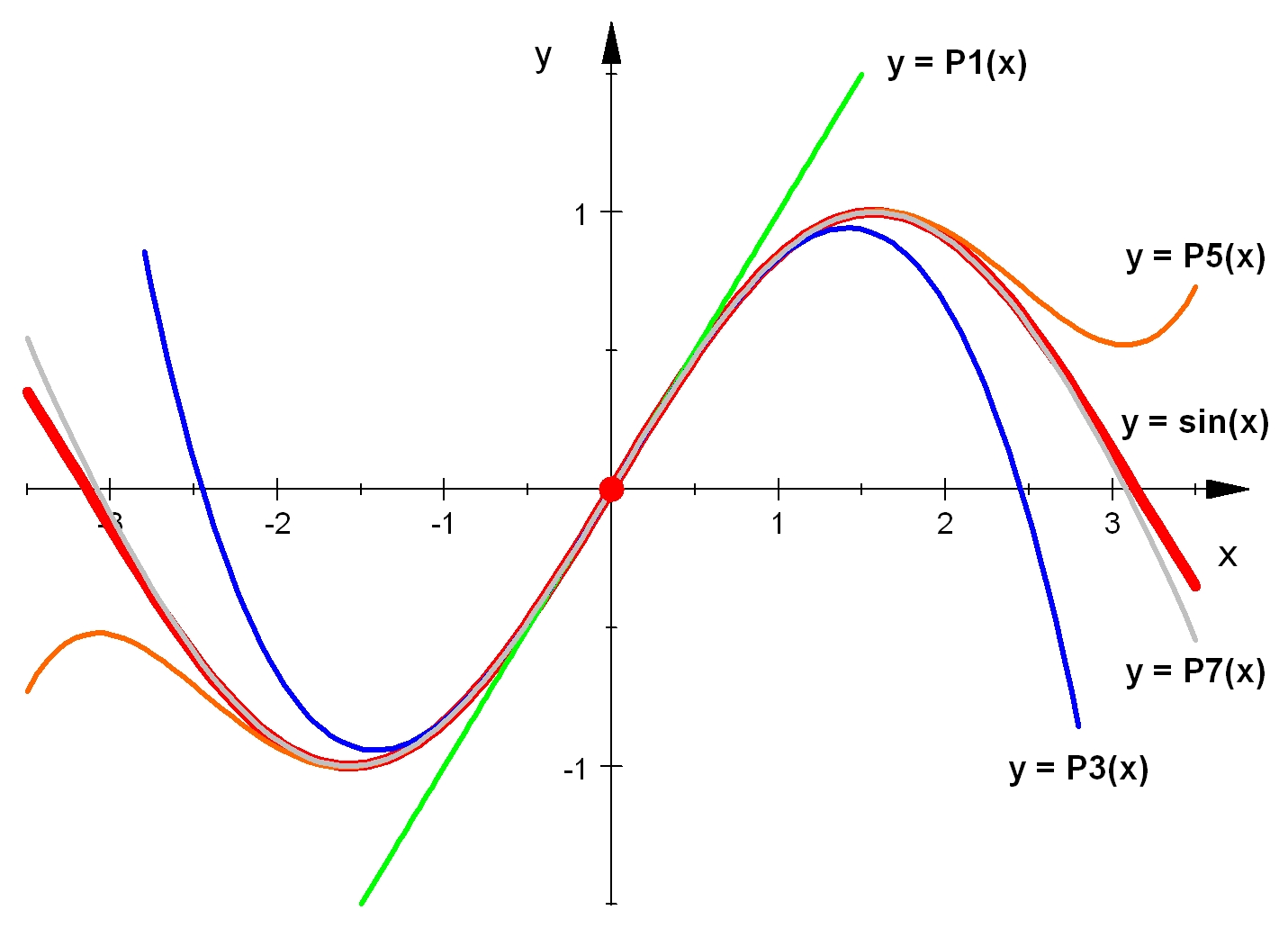}
\caption{Taylor Approximations to $y=\sin x$}        
\label{fig910}
 \end{figure}

 \section{Conclusion}  
 Adaptive Horizon Model Predictive Control is a scheme for varying the horizon length  in Model Predictive Control as the stabilization process evolves. It adapts  the horizon in real time by testing Lyapunov and feasibility conditions on  extensions of  optimal trajectories returned by the nonlinear program solver.   In this way it seeks the shortest horizons consistent with stabilization.
 
 AHMPC requires a terminal cost and terminal feedback that stabilizes the plant in some neighborhood of the operating point but that neighborhood need not be known explictly.   Higher degree Taylor polynomial approximations to the optimal cost and the optimal feedback of the coresponding infinite horizon optimal control problems can be found by an extension of Al'brekht's method \cite{AK12}.     The higher degree Taylor polynomial approximations to optimal cost need not be positive definite but they can be extended to  nonnegative definite polynomials by completing the squares.   These nonnegative definite extensions  and the Taylor polynomial approximations to the optimal feedback can be used as terminal costs and terminal feedbacks in AHMPC.   We have shown by an example  that a higher degree terminal cost and feedback can outperform using LQR to define a degree two terminal cost and a degree one terminal feedback.


\end{document}